# Symmetry-breaking instabilities of convection in squares

By A.M. Rucklidge

*Department of Applied Mathematics and Theoretical Physics,
University of Cambridge, Cambridge CB3 9EW, UK*

Convection in an infinite fluid layer is often modelled by considering a finite box with periodic boundary conditions in the two horizontal directions. The translational invariance of the problem implies that any solution can be translated horizontally by an arbitrary amount. Some solutions travel, but those solutions that are invariant under reflections in both horizontal directions cannot travel, since motion in any horizontal direction is balanced by an equal and opposite motion elsewhere. Equivariant bifurcation theory allows us to understand the steady and time-dependent ways in which a pattern can travel when a mirror symmetry of the pattern is broken in a bifurcation. Here we study symmetry-breaking instabilities of convection with a square planform. A pitchfork bifurcation leads to squares that travel uniformly, while a Hopf bifurcation leads to a new class of oscillations in which squares drift to and fro but with no net motion of the pattern. Two types of travelling squares are possible after a pitchfork bifurcation, and three or more oscillatory solutions are created in a Hopf bifurcation. One of the three oscillations, alternating pulsating waves, has been observed in recent numerical simulations of convection in the presence of a magnetic field. We also present a low-order model of three-dimensional compressible convection that contains these symmetry-breaking instabilities. Our analysis clarifies the relationship between several types of time-dependent patterns that have been observed in numerical simulations of convection.

## 1. Introduction

Symmetries play a central role in the analysis of the dynamics of convection. Often an infinite layer is modelled by a periodic box, and the range of allowed behaviour is influenced by the symmetries of the box. This approach is justified by the experimental observation of convection in periodic patterns: rolls, squares (Le Gal & Croquette 1988) and hexagons. The presence of the underlying translational invariance of the layer implies that periodic spatial patterns are not isolated fixed points of the dynamical system, but that there are continuous families (or group orbits) of fixed points generated by translations. When the spatial pattern of convection is invariant under a reflection symmetry, the pattern cannot drift, but symmetry-breaking bifurcations lead to patterns that drift slowly along their group orbit.

A two-dimensional periodic box has the symmetry $O(2) = Z_2 \ltimes SO(2)$, where $\ltimes$ denotes the semi-direct product; two-dimensional convection initially takes the form of $Z_2$-symmetric rolls, which are invariant under reflection in the vertical plane between pairs of counter-rotating rolls. The $SO(2)$ translational invariance of the periodic box implies that there is a continous group orbit of rolls parameterised by their position in the box. A pitchfork bifurcation that breaks the mirror symmetry of the rolls leads to patterns that travel steadily in one direction along the group orbit, while a symmetry-breaking Hopf bifurcation leads to patterns that travel to and fro with no net motion (Landsberg & Knobloch 1991; Proctor & Weiss 1993; Matthews *et al.* 1993). These solutions are known as travelling waves and pulsating waves respectively. In numerical simulations, the instability that breaks the mirror symmetry is associated with the generation of a horizontal shear across the layer.

We consider three-dimensional convection in a plane layer of unit depth with periodic boundary conditions with equal wave numbers $k$ in the horizontal directions $x$ and $y$, in the presence of an imposed vertical magnetic field. The system has the symmetry group $\Gamma = D_4 \ltimes T^2$, where $D_4$ is the symmetry group of a square and $T^2$ is the two-torus of translations in the two horizontal directions. There is a trivial solution that is invariant under $\Gamma$; in the prescribed periodic box, the trivial solution loses stability to rolls and squares, which break the translation invariance but preserve some of the symmetries in $D_4$. There are two circles of $Z_2 \times O(2)$-invariant rolls, related to each other by a reflection in the diagonal, and there is a two-torus of $D_4$-invariant squares, parameterised by their position in the box.

We confine ourselves to spatially periodic behaviour with the same period as the basic pattern; long-wavelength instabilities of squares have already been examined (Hoyle 1993). Our work is complementary to that of McKenzie (1988), who considered instabilities of square and hexagonal patterns that preserve a mirror symmetry and thus do not lead to travelling patterns. McKenzie also considered steady spatial period-doubling transitions.

In section two we discuss symmetry-breaking pitchfork and Hopf bifurcations from squares. In section three, we present a model of three-dimensional convection; in this model, it is a spontaneously generated shear across the convecting layer that is responsible for driving the underlying square pattern. We discuss the relevance of our analysis to other problems in the final section.

## 2. Normal forms for secondary instabilities of squares

The presentation will be couched in the language of convection, but the analysis in this section is applicable more generally. The state of the convecting system is specified by the velocities, temperature, density and magnetic field; these are collected into a single vector $U(x, y, z, t)$. The dynamics of the system is governed by a set of partial differential equations (PDEs) that we write symbolically as

$$\frac{dU}{dt} = \mathcal{F}(U; \lambda). \tag{2.1}$$

The parameter $\lambda$ represents the Rayleigh number, proportional to the temperature difference across the layer.

The full symmetry group of the system is $\Gamma = D_4 \ltimes T^2$. The group $D_4$ (the point symmetries of a square) is generated by reflections in two vertical planes: $m_x$, which sends $x$ to $-x$, and $m_d$, which exchanges $x$ and $y$. We define $r_q$, a rotation by one quarter turn, by $r_q = m_d m_x$, $m_y$, the reflection that sends $y$ to $-y$, by $m_y = m_d m_x m_d$, and $m_{d'}$, the reflection in the other diagonal, by $m_{d'} = m_x m_d m_x$. The translation $\tau(\delta_x, \delta_y)$, an element of the two-torus $T^2$, shifts $x$ and $y$ by $\delta_x$ and $\delta_y$ respectively. The dynamics specified by $\mathcal{F}$ is equivariant under the action of these

symmetries, so
$$\mathcal{F}(\gamma U; \lambda) = \gamma \mathcal{F}(U; \lambda), \tag{2.2}$$
where $\gamma \in \Gamma$. The symmetry operation $\gamma$ acts on a given solution $U$ to generate a new (or perhaps the same) solution.

We consider a steady solution $U_0(x, y, z)$ that describes convection with a square planform. This solution is invariant under the $D_4$ spatial group, so
$$m_x U_0 = m_d U_0 = U_0, \tag{2.3}$$
and the translations $\tau(\delta_x, \delta_y)$ generate a continuous family of square solutions parameterised by the shifts $\delta_x$ and $\delta_y$. This two-torus of solutions is known as the group orbit through $U_0$. Since $U_0$ is an equilibrium point, as are all states on the group orbit,
$$\mathcal{F}(U_0; \lambda) = 0 = \mathcal{F}(\tau(\delta_x, \delta_y)U_0; \lambda) \tag{2.4}$$
for any $(\delta_x, \delta_y) \in [0, 2\pi/k)^2$.

The eigenvalues of the Jacobian matrix evaluated at $U_0$ determine the stability of the square equilibrium point. Because the group orbit is continuous, the square pattern is neutrally stable with respect to translations in the $x$ and $y$ directions; this implies that the Jacobian matrix always has a pair of zero eigenvalues corresponding to the two horizontal directions. We assume that the square pattern undergoes a simple symmetry-breaking instability at $\lambda = 0$, at which point there will be additional zero or pure imaginary eigenvalues. By simple, we mean that the eigenvalues that cross the imaginary axis are no more degenerate than required by the symmetry of the problem. By symmetry-breaking, we mean that the bifurcating solutions have less symmetry than the original square pattern. We do not consider instabilities that break the initial spatial periodicity.

In order to study the behaviour near the instability, we expand the solution $U$ near the square equilibrium $U_0$:
$$U(x, y, z, t) = \tau(\phi_x(t), \phi_y(t))\{U_0(x, y, z) + V(x, y, z, t)\}. \tag{2.5}$$

Here we are following an approach used to study instabilities of vortices in the Taylor–Couette problem (Iooss 1986). The idea is that we require an expansion that is valid in the neighbourhood of the full group orbit of $U_0$; when $V$ is small in (2.5), the system will be near a steady solution that is a translation of the original square pattern. The two variables $\phi_x$ and $\phi_y$ are scalar functions of time that represent the phase of the square pattern in the periodic domain; these are not required to be small. The function $V$ is a small vector function of space and time that describes how the system differs from squares, but $V$ does not contain any translations. With this expansion, the symmetries act in the following way:
$$m_x(\phi_x, \phi_y, V) = (-\phi_x, \phi_y, m_x V), \tag{2.6}$$
$$m_d(\phi_x, \phi_y, V) = (\phi_y, \phi_x, m_d V), \tag{2.7}$$
$$\tau(\delta_x, \delta_y)(\phi_x, \phi_y, V) = (\phi_x + \delta_x, \phi_y + \delta_y, V). \tag{2.8}$$

Substituting the expansion (2.5) into the PDE (2.1), we recover two ordinary differential equations (ODEs) for the translation variables $\phi_x$ and $\phi_y$ and a PDE for $V$:
$$\dot{\phi}_x = f_x(V), \qquad \dot{\phi}_y = f_y(V), \qquad \frac{dV}{dt} = \tilde{\mathcal{F}}(V; \lambda), \tag{2.9}$$

where $f_x$ and $f_y$ are real functions and $\tilde{\mathcal{F}}$ represents the projection of the original PDE onto the space with all translations factored out. The dot indicates differentiation with respect to time in the two ODEs. The variables $\phi_x$ and $\phi_y$ do not appear on the right-hand sides of (2.9) since the dynamics must be invariant under the action of translations. The mirror symmetries $m_x$ and $m_d$ impose restrictions on the form of (2.9):

$$f_x(m_x V) = -f_x(V), \qquad f_y(m_x V) = f_y(V), \qquad f_x(m_d V) = f_y(V). \tag{2.10}$$

In addition, $\tilde{\mathcal{F}}$ must be $D_4$-equivariant. Thus we have reduced the PDE for $U$, with $D_4 \ltimes T^2$ symmetry, to a PDE for $V$ with $D_4$ symmetry along with a pair of ODEs that describe how the square pattern drifts. The $T^2$ torus of translations now acts in a simple way on the variables $\phi_x$ and $\phi_y$.

The original square pattern, with $V = 0$ and $\phi_x$ and $\phi_y$ constant, is neutrally stable with respect to translations; the two zero eigenvalues forced by this appear in the two ODEs in (2.9). Squares lose stability when eigenvalues of the Jacobian matrix of $\tilde{\mathcal{F}}$ cross the imaginary axis. We assume that there is a simple symmetry-breaking bifurcation at $\lambda = 0$. This could be either a pitchfork bifurcation (zero eigenvalues) or a Hopf bifurcation (pure imaginary eigenvalues). These eigenvalues will have multiplicity two because of the $D_4$ symmetry.

At the bifurcation point, a centre manifold reduction of the reduced PDEs yields the normal forms for the pitchfork and the Hopf bifurcations with $D_4$ symmetry, augmented by two ODEs for $\phi_x$ and $\phi_y$ (Chossat & Iooss 1994). In both cases, let $v_x(t)$ and $v_y(t)$ be the amplitudes of the critical modes; on the centre manifold, $V$ will be a function of these amplitudes. The amplitudes $v_x$ and $v_y$ will be real in the steady case and complex in the oscillatory case. The translations act trivially on $v_x$ and $v_y$. The modes can be chosen so that the mirror symmetries act on $v_x$ and $v_y$ in the following way:

$$m_x(v_x, v_y) = (-v_x, v_y), \qquad m_d(v_x, v_y) = (v_y, v_x). \tag{2.11}$$

Together with the action of the reflections on $\phi_x$ and $\phi_y$, this further restricts the form of the functions $f_x$ and $f_y$ in the ODEs:

$$f_x(v_x, v_y) = \operatorname{Re}(f(v_x^2, v_y^2) v_x), \qquad f_y(v_x, v_y) = \operatorname{Re}(f(v_y^2, v_x^2) v_y), \tag{2.12}$$

with $f$ a smooth function that may be complex in the case of a Hopf bifurcation. The form of the functions $f_x$ and $f_y$ implies that the original squares will drift in the $x$ direction when $v_x \neq 0$ and in the $y$ direction when $v_y \neq 0$.

Thus when squares lose stability in a pitchfork bifurcation at $\lambda = 0$, the behaviour near the bifurcation is governed by

$$\dot{v}_x = \lambda v_x + A v_x^3 + B v_y^2 v_x, \qquad \dot{\phi}_x = D v_x, \tag{2.13}$$

$$\dot{v}_y = \lambda v_y + A v_y^3 + B v_x^2 v_y, \qquad \dot{\phi}_y = D v_y, \tag{2.14}$$

which is the normal form for a pitchfork bifurcation with $D_4$ symmetry, augmented by two equations for the drift along the group orbit of the original square solution. The real constants $A$, $B$ and $D = f(0,0)$ will depend on $\lambda$ and $U_0$, and can in principle be computed from the original PDEs.

There are two types equilibria created in this pitchfork bifurcation. Travelling squares (TSq) travel in, say, the $x$ direction, with $v_y = \dot{\phi}_y = 0$, $v_x \neq 0$ and $\dot{\phi}_x =$ constant; these are invariant under $m_y$ and drift uniformly along the group orbit of squares induced by $x$ translations. Diagonally travelling squares (DTSq) have $v_x = v_y$ and $\dot\phi_x = \dot\phi_y =$ constant; these are invariant under $m_d$ and drift diagonally. The properties of these travelling solutions are summarised below (isotropy subgroups, up to a conjugacy, are also indicated):

$$\text{TSq:} \quad v_x \neq 0, \quad v_y = 0, \quad \langle m_y \rangle,$$
$$\text{DTSq:} \quad v_x = v_y \neq 0, \quad \langle m_d \rangle,$$

where angle brackets denote the subgroup of $D_4$ generated by the given elements, isomorphic to $Z_2$ in each case. Conjugate solutions are generated by the $D_4$ symmetries. The stability of these solutions is determined by the parameters $A$ and $B$.

The other possibility is that squares lose stability in a Hopf bifurcation with frequency $\omega$ at $\lambda = 0$. Here the behaviour is governed by

$$\dot v_x = \left((\lambda + \mathrm{i}\omega) + A(|v_x|^2 + |v_y|^2) + B|v_x|^2\right) v_x + C\bar v_x v_y^2, \qquad (2.15)$$
$$\dot v_y = \left((\lambda + \mathrm{i}\omega) + A(|v_x|^2 + |v_y|^2) + B|v_y|^2\right) v_y + C\bar v_y v_x^2, \qquad (2.16)$$
$$\dot\phi_x = \mathrm{Re}(Dv_x), \qquad (2.17)$$
$$\dot\phi_y = \mathrm{Re}(Dv_y), \qquad (2.18)$$

which is the normal form for a Hopf bifurcation with $D_4$ symmetry (Swift 1988), augmented by two equations for the drift of the pattern. Here $A$, $B$, $C$ and $D$ are complex constants that will again depend on $\lambda$ and $U_0$.

This normal form is $D_4 \ltimes T^2$ equivariant, but it is also equivariant under increasing the phases of the complex amplitudes $v_x$ and $v_y$ together by the same amount. This additional $S^1$ phase symmetry of the normal form is not a symmetry of the full convection problem: solutions of the Hopf normal form (2.15–2.18) make up a continuous family generated by multiplying $v_x$ and $v_y$ by complex numbers on the unit circle, while solutions of the convection problem do not share this property. The additional symmetry arises because non-$S^1$ equivariant nonlinear terms can be eliminated from the normal form by near-identity transformations (see Swift (1988) for more details).

However, some elements of $S^1$, which we interpret as discrete temporal phase shifts, are important in describing the symmetry properties of periodic orbits in the full convection problem. Consider the effect of a symmetry operation $\gamma \in D_4$ on a periodic orbit $\mathbf{v}(t)$ with period $T$, where $\mathbf{v} = (v_x, v_y)$. Suppose $\gamma^n$ is the identity ($n = 2$ for reflections and $n = 4$ for rotations by $90°$). There are two possible outcomes: $\mathbf{v}(t)$ and $\gamma \mathbf{v}(t)$ are either disjoint or identical periodic orbits, with possibly a phase shift in time (Golubitsky *et al.* 1988). In the latter case, we must have $\gamma \mathbf{v}(t) = \mathbf{v}(t + \frac{m}{n}T)$, where $m$ is an integer; then $\gamma^n \mathbf{v}(t) = \mathbf{v}(t + mT) = \mathbf{v}(t)$. Therefore we define two symmetry operations $t_h, t_q \in S^1$, which advance time by one half or one quarter of the period, respectively. The other elements of $S^1$ do not carry over to the full problem, but the discrete time shifts form part of the spatio-temporal symmetries that describe the periodic orbits created in the Hopf bifurcation in the full problem.

The $D_4$-Hopf normal form (2.15–2.16) without the auxiliary drift equations has been studied in the context of bifurcations from the trivial solution (Swift 1988). Three periodic orbits are always created in this Hopf bifurcation, and there is the additional possibility of two other types of solution, one periodic and the second doubly periodic, for some parameter values. Demonstrating the existence of the doubly periodic torus requires the $S^1$ symmetry of the Hopf normal form, so the torus may not persist in the full problem. In the case of a $D_4$-Hopf bifurcation from the trivial solution, the three types of periodic orbit are named standing rolls, standing squares and alternating rolls (Silber & Knobloch 1991). The fourth periodic solution, which does not exist for all parameter values, is called standing cross rolls, and is always unstable near the bifurcation point. In our context, where the Hopf bifurcation is from a nontrivial solution with $D_4$ symmetry, we name the periodic orbits pulsating squares (PSq), diagonally pulsating squares (DPSq) and alternating pulsating waves (APW). The names are chosen to describe solutions that have been observed in numerical simulations of compressible magnetoconvection (Matthews *et al.* 1994; Matthews *et al.* 1995). The periodic orbits are characterised by the following relations between their amplitudes and isotropy subgroups:

$$\begin{aligned}\text{PSq:} \quad & v_x \neq 0, \quad v_y = 0, \quad & \langle m_y, t_h m_x \rangle, \\ \text{DPSq:} \quad & v_x = v_y \neq 0, \quad & \langle m_d, t_h m_{d'} \rangle, \\ \text{APW:} \quad & v_x = \mathrm{i} v_y \neq 0, \quad & \langle t_q r_q \rangle,\end{aligned}$$

where angle brackets denote the subgroup of $D_4 \times S^1$ generated by the given elements. The isotropy subgroups in each case have four elements and are isomorphic to $D_2$, $D_2$ and $Z_4$, respectively. The stability of the three periodic orbits as the parameters in the normal form are varied has been computed (Silber & Knobloch 1991); the orbits may bifurcate subcritically or supercritically, and none, one or two may be stable. In particular, all three may bifurcate supercritically with all three unstable.

Since these periodic orbits all satisfy $t_h \mathbf{v} = -\mathbf{v}$, the translation variables $\phi_x$ and $\phi_y$ will be driven in one direction in the first half of the oscillation and in the other direction in the second half, resulting in the underlying square pattern drifting to and fro but returning exactly to its original location after one period. This oscillating drift is along the coordinate axes for the pulsating squares and along the diagonals for the diagonal pulsating squares. In the case of alternating pulsating waves, the pattern drifts first in the positive $x$ direction (say), then in the positive $y$ direction, then in the negative $x$ direction, and finally back in the negative $y$ direction; there is another alternating pulsating wave that circulates in the other direction. In all three cases, in common with the situation in two-dimensional convection (Landsberg & Knobloch 1991; Proctor & Weiss 1993), there is no net drift over the course of the oscillation. Using the symmetries of the other two oscillatory solutions, we can also deduce that there will be no net drift in the case of the oscillation corresponding to standing cross rolls (though the drift will have different amplitudes in the two horizontal directions), and that on average, there will be no net drift in the case of torus.

## 3. Example: a model of compressible magnetoconvection

We present as an example a set of ODEs that is equivariant under $\Gamma$ and that exhibits the bifurcations we have been discussing. The ODEs are based on a model of three-dimensional incompressible convection in an imposed vertical magnetic field (Rucklidge & Matthews 1995), extended to include the effects of compressibility by including terms suggested by work on the corresponding two-dimensional problem (Proctor *et al.* 1994). The modes used in the model are shown in figure 1. The variables $a_x$ and $a_y$ are complex and represent the amplitudes and positions of the $x$- and $y$-rolls in the periodic box (figure 1a,b); squares are represented by an equal sum of $x$- and $y$-rolls. The real variables $c_x$ and $c_y$ represent shearing modes: the horizontal averages of the two components of the velocity are in opposite directions at the top and bottom of the layer (figure 1c,d). Rolls and squares formed in the initial bifurcation have $c_x = c_y = 0$. We also include real variables $d_x$ and $d_y$ to represent the $x$ and $y$ components of a magnetic field imposed across the layer; the magnetic field is stretched out horizontally by the shear. Other modes, with horizontal dependence, are also involved in these symmetry-breaking bifurcations, but the analysis in terms of horizontally averaged modes is more transparent and more natural since the shearing modes decay the slowest.

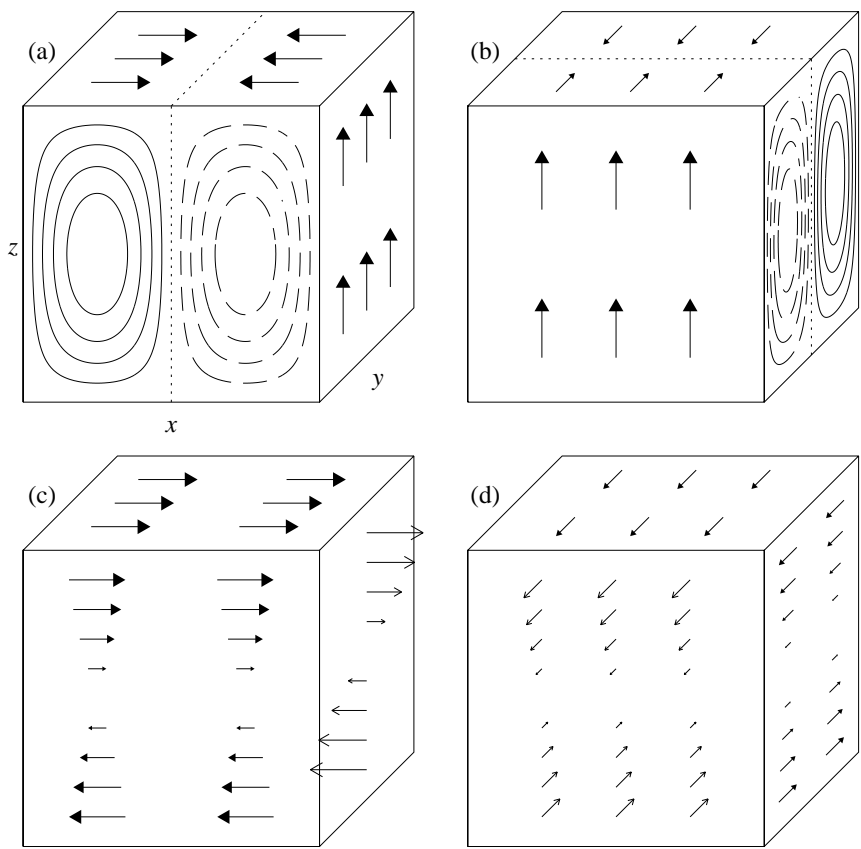

Figure 1. The four modes used to represent convection: (a) $x$-rolls ($a_x$); (b) $y$-rolls ($a_y$); (c) $x$-shear ($c_x$); (d) $y$-shear ($c_y$). Arrows with solid heads represent fluid motion parallel to the sides or tops of the cubes, and arrows with open heads indicate fluid motion that crosses the boundary of the cube. In (a) and (b), the solid lines represent clockwise circulations and the broken lines counter-clockwise circulations, while the dotted lines show the vertical mirror planes of symmetry.

The action of the elements of $\Gamma$ on the mode amplitudes is given by

$$m_x(a_x, c_x, d_x, a_y, c_y, d_y) = (\bar{a}_x, -c_x, -d_x, a_y, c_y, d_y), \tag{3.1}$$

$$m_d(a_x, c_x, d_x, a_y, c_y, d_y) = (a_y, c_y, d_y, a_x, c_x, d_x), \tag{3.2}$$

$$\tau(\delta_x, \delta_y)(a_x, c_x, d_x, a_y, c_y, d_y) = (a_x e^{ik\delta_x}, c_x, d_x, a_y e^{ik\delta_y}, c_y, d_y). \tag{3.3}$$

Thus the $m_x$ reflection changes the sign of the phase of $a_x$, the direction of the $x$ shear and the direction of the $x$ component of the magnetic field; the $m_d$ reflection exchanges $x$ and $y$; and translations change the phases of $a_x$ and $a_y$.

In the convection problem, the trivial conducting state loses stability as the controlling parameter, the Rayleigh number $R$, increases above a critical value $R_C$. Above this critical value, convection takes the form of rolls or squares. This behaviour is represented by the normal form for the $\Gamma$-symmetric pitchfork bifurcation:

$$\dot{a}_x = \mu a_x - |a_x|^2 a_x - (1+\beta)|a_y|^2 a_x, \tag{3.4}$$
$$\dot{a}_y = \mu a_y - |a_y|^2 a_y - (1+\beta)|a_x|^2 a_y. \tag{3.5}$$

Here, $\mu$ and $\beta$ are real parameters, with $\mu$ proportional to $R - R_C$. Rolls ($|a_x|^2 = \mu$, $a_y = 0$) and squares ($|a_x|^2 = |a_y|^2 = \mu/(2+\beta)$) are created in the pitchfork bifurcation at $\mu = 0$, which is assumed to be supercritical. Rolls are stable if $\beta > 0$ while squares are stable if $-2 < \beta < 0$. The two-torus of $D_4$-invariant squares is parameterised by the phases of $a_x$ and $a_y$.

Rolls and squares subsequently lose stability to a symmetry-breaking shearing mode, which is modelled by including $c_x$, $d_x$, $c_y$ and $d_y$ in the above ODEs:

$$\dot{a}_x = \mu a_x - |a_x|^2 a_x - (1+\beta)|a_y|^2 a_x - \gamma c_x^2 a_x + iDa_x c_x, \tag{3.6}$$
$$\dot{c}_x = -c_x - Qd_x + |a_x|^2 c_x, \tag{3.7}$$
$$\dot{d}_x = \zeta c_x - \zeta d_x, \tag{3.8}$$
$$\dot{a}_y = \mu a_y - |a_y|^2 a_y - (1+\beta)|a_x|^2 a_y - \gamma c_y^2 a_y + iDa_y c_y, \tag{3.9}$$
$$\dot{c}_y = -c_y - Qd_y + |a_y|^2 c_y, \tag{3.10}$$
$$\dot{d}_y = \zeta c_y - \zeta d_y. \tag{3.11}$$

These ODEs are equivariant under the action of $\Gamma$ given by (3.1–3.3). All the parameters are required by the symmetry to be real. The strength of the imposed magnetic field, which acts to oppose the shear, is given by $Q$, $\zeta$ is the Ohmic diffusivity divided by the viscosity, $\gamma$ describes how the shear damps the vigour of convection, and $D$ determines how the rolls and squares will drift once the $D_4$ symmetry is broken.

The phases of $a_x$ and $a_y$ decouple and the model can be rewritten in the form

$$\dot{r}_x = \mu r_x - r_x^3 - (1+\beta)r_y^2 r_x - \gamma c_x^2 r_x, \qquad \dot{\theta}_x = Dc_x, \tag{3.12}$$
$$\dot{c}_x = -c_x - Qd_x + r_x^2 c_x, \tag{3.13}$$
$$\dot{d}_x = \zeta c_x - \zeta d_x, \tag{3.14}$$
$$\dot{r}_y = \mu r_y - r_y^3 - (1+\beta)r_x^2 r_y - \gamma c_y^2 r_x, \qquad \dot{\theta}_y = Dc_y, \tag{3.15}$$
$$\dot{c}_y = -c_y - Qd_y + r_y^2 c_y, \tag{3.16}$$
$$\dot{d}_y = \zeta c_y - \zeta d_y, \tag{3.17}$$

where $a_x = r_x e^{i\theta_x}$ and $a_y = r_y e^{i\theta_y}$. Comparing the equations for the phases of $a_x$ and $a_y$ above with the equations for the displacements along the group orbit in the normal forms (2.13–2.14) and (2.17–2.18), it is clear how the shearing mode drives the drift of the square pattern: the rate of drift in the $x$ direction ($\dot{\theta}_x$) is proportional to the shear amplitude $c_x$. The connection between this model and the general theory in the previous section is also apparent. The underlying square pattern ($U_0$) is described by $r_x$ and $r_y$; the perturbation from squares $(V)$ is given by $c_x$, $c_y$, $d_x$ and $d_y$, and these perturbations drive the square pattern along its group orbit in the same way as in (2.9).

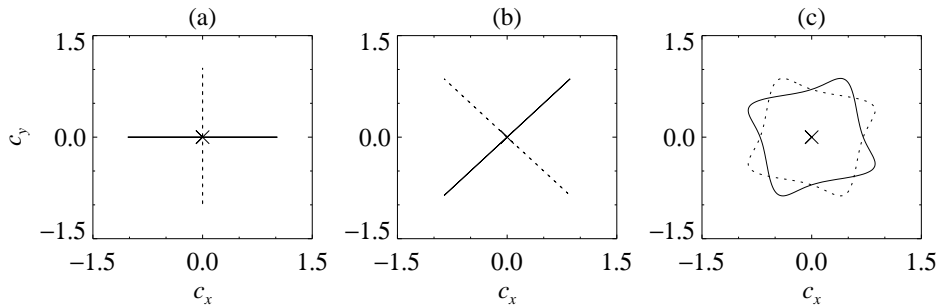

Figure 2. Examples of the three periodic orbits created in the Hopf bifurcation in (3.6–3.11) with $\mu = 1.0$, $Q = 1.0$, $\zeta = 0.2$, $\beta = -1.5$ and $\gamma = 1.0$. For these parameter values, the Hopf bifurcation occurs at $\mu = 0.6$. Also shown (dotted lines) are conjugate orbits. (a) Pulsating squares (PSq), invariant under $m_y$ and $t_h m_x$; (b) Diagonally pulsating squares (DPSq), invariant under $m_d$ and $t_h m_{d'}$; (c) Alternating pulsating waves (APW), invariant under $t_q r_q$. The first two orbits are unstable and the last one (APW) is stable. The cross represents the original square fixed point, with $c_x = c_y = 0$.

In the model, squares undergo a pitchfork bifurcation when $\mu = (1+Q)(2+\beta)$ and a Hopf bifurcation when $\mu = (1+\zeta)(2+\beta)$ with $Q > \zeta$. The pitchfork bifurcation leads to uniformly drifting TSq and DTSq, while the Hopf bifurcation leads to time-dependent patterns that drift to and fro with no net drift. With illustrative parameter values, two of the three periodic orbits created in this Hopf bifurcation are unstable and one (APW) is stable (see figure 2). Other possibilities exist for other parameter values; an exhaustive study of the behaviour in different parameter regimes is beyond the scope of this paper. Note how the APW (figure 2c) are not circles (and so not invariant under the $S^1$ symmetry of the Hopf normal form), but are unchanged under $t_q r_q$. There is also a codimension-two Takens–Bogdanov bifurcation with $D_4$ symmetry when $(\mu, Q) = ((1+\zeta)(2+\beta), \zeta)$, at which point the pitchfork and Hopf bifurcations coincide.

## 4. Discussion

We have considered the symmetry-breaking bifurcations from convection in a square planform. The overall symmetry group is $\Gamma = D_4 \ltimes T^2$, but the $D_4$-invariant squares break the $T^2$ symmetry, so the symmetry group relevant to the secondary bifurcations is just $D_4$. However, the presence of the $T^2$ symmetry implies that the bifurcating solutions will drift along the group orbit of the underlying squares. Thus in a pitchfork bifurcation, the squares will travel along the coordinates axes or diagonally, and in a Hopf bifurcation, squares will travel to and fro, along the coordinate axes, diagonally or alternating between the two coordinate directions. There is the possibility of additional bifurcating solutions in the oscillatory case: an additional periodic solution with a smaller symmetry group (equivalent to the standing cross rolls of the bifurcation from the trivial solution) and a doubly periodic solution, which may be attracting but whose existence relies on the $S^1$ symmetry introduced by the Hopf normal form (2.15–2.18). The full convection problem does not have this $S^1$ symmetry (nor do the model ODEs (3.6–3.11)), so it is likely that if a torus is created in the Hopf bifurcation then there will be phase locking far enough away from the bifurcation. There will be further complications (for example, homoclinic tangencies) if the torus is involved in global bifurcations (Swift 1988).

Boussinesq convection with the same boundary condition at the top and bottom of the layer differs from compressible convection in that there is an additional mirror symmetry of reflections in the horizontal midplane. This implies that there will be two types of vertical shear eigenfunction, distinguished by being even (velocities in the same direction at the top and bottom of the layer) or odd (velocities in the opposite direction) under reflection in the horizontal midplane. Both parities have been used in models of two-dimensional magnetoconvection (Matthews *et al.* 1993; Landsberg & Knobloch 1993), and the odd parity has been used in the model of three-dimensional magnetoconvection (Rucklidge & Matthews 1995) that forms the basis of the ODE model presented in section 3. The interactions between the roll and shear modes in the two cases are different: in the even case, rolls and shear interact directly, but in the odd case, they interact through a third mode that represents the tilting of the original rolls or squares. The analysis in this paper applies to the even and odd cases separately, but there are complications if both even and odd modes are excited. In addition, in the odd case, velocities at the top and bottom of the layer will be equal and opposite, implying that the patterns will not travel; the coefficient $D$ in (2.13–2.14) and (2.17–2.18) must therefore be zero in this case.

Steady squares are the preferred planform at the onset of compressible convection in a strong magnetic field with Ohmic diffusivity greater than the thermal diffusivity, and some shearing instabilities of squares have been observed but not explored in detail (Matthews *et al.* 1995); our analysis is likely to be relevant in this case. The transition from squares to travelling squares (TSq and DTSq) has been observed in a study of $\sqrt{2} : 1$ resonance in non-Boussinesq convection (Proctor & Matthews 1996). In addition, clear examples of alternating pulsating waves have been reported in numerical simulations of convection in a weak magnetic field (Matthews *et al.* 1994; Matthews *et al.* 1995), though these were created not in a Hopf bifurcation but after a complicated sequence of local and global bifurcations.

In the two-dimensional problem, pulsating waves (rolls that travel to and fro, with a spatio-temporal symmetry $t_h m_x$) can be created in a symmetry-breaking Hopf bifurcation from rolls, or after a sequence of bifurcations. Rolls can undergo a symmetry-breaking pitchfork bifurcation to travelling rolls, which subsequently lose stability to modulated travelling rolls in a Hopf bifurcation; finally, the $t_h m_x$ symmetry is restored when pulsating waves are created in a global bifurcation. These two routes to pulsating waves are brought together at a codimension-two Takens–Bogdanov bifurcation point, where the Hopf and pitchfork bifurcations coincide (Matthews *et al.* 1993). The use of normal forms and low-order models has proved essential in this two-dimensional analysis (Rucklidge & Matthews 1996).

The situation in three dimensions is more complicated for two reasons. First, near the Takens–Bogdanov point, there are two types of travelling squares and at least three types of oscillating squares to be fitted into a coherent picture. A preliminary investigation has revealed that there is also chaotic dynamics near the codimension-two point (Armbruster *et al.* 1989), and a detailed study will be necessary to elucidate the many routes from travelling solutions to solutions that oscillate to and fro. A study of the Takens–Bogdanov bifurcation with $D_4$ symmetry is also of interest as it is an important subcase of Takens–Bogdanov with $D_4 \ltimes T^2$, which arises in the context of the initial instability in magnetoconvection at parameter values where the pitchfork bifurcation to steady convection and the Hopf bifurcation to oscillatory convection coincide (Clune & Knobloch 1994; Matthews *et al.* 1995).

The second complication in three dimensions is that alternating pulsating waves have been observed in numerical simulations of convection in the presence of a vertical magnetic field at parameter values for which rolls, not squares, are preferred at onset (Matthews *et al.* 1994; Matthews *et al.* 1995). The route from rolls to APW is a very involved sequence of local and global bifurcations that depends in detail on the parameter values. For example, in the numerical simulations of magnetoconvection, steady rolls lose their mirror symmetry, become time dependent then three dimensional, and finally gain the $t_q r_q$ symmetry in a global bifurcation to form APW (Matthews *et al.* 1995). A similar sequence (but including structurally stable heteroclinic cycles) also occurs in (3.6–3.11) and in the model of Boussinesq magnetoconvection on which those equations are based (Rucklidge & Matthews 1995), and a sequence of bifurcations leading from rolls to modulated diagonally travelling squares has been described in simulations of non-magnetic convection (Matthews *et al.* 1996). In all cases, low-order models have clarified the global bifurcations and the changes in symmetry that occur in the numerical simulations.

It is a pleasure to acknowledge the discussions I have had with Mary Silber, which greatly increased my understanding of this subject. I have also benefited from discussions with Rebecca Hoyle, Edgar Knobloch, Paul Matthews, Michael Proctor, and Nigel Weiss. This research was supported by SERC and its successor EPSRC, and by Peterhouse, Cambridge.